\documentclass[8pt]{article}
\usepackage{mathrsfs}
\usepackage{amsthm}
\usepackage{amssymb}
\usepackage{amsmath}
\usepackage{graphicx}
\usepackage{color}
\usepackage{amsfonts}
\usepackage{float}
\usepackage{cite}
\usepackage [latin1]{inputenc}
\usepackage[text={140mm,210mm},left=35mm,vmarginratio=1:1]{geometry}
\newtheorem{theorem}{Theorem}[section]

\newtheorem{lemma}[theorem]{Lemma}
\newtheorem{proposition}[theorem]{Proposition}

\numberwithin{equation}{section}
\normalsize

\begin{document}
\title{\textbf{Nonequilibrium moderate deviations from hydrodynamics of simple exclusion processes}}

\author{Xiaofeng Xue \thanks{\textbf{E-mail}: xfxue@bjtu.edu.cn \textbf{Address}: School of Mathematics and Statistics, Beijing Jiaotong University, Beijing 100044, China.}\\ Beijing Jiaotong University}

\date{}
\maketitle

\noindent {\bf Abstract:} In this paper, we give the moderate deviation principle from the hydrodynamic limit of the simple exclusion process on $1$-dimensional torus starting from a nonequilibrium state, which extends the result given in Gao and  Quastel (2003) about the case where the process starts from an equilibrium state. The exponential tightness of the scaled density field of the process and a replacement lemma play key roles in the proof of the main result. We utilize Grownwall's inequality and the upper bound of the large deviation principle given in Kipnis, Olla and Varadhan (1989) to prove above exponential tightness and replacement lemma respectively in the absence of the invariance of the initial distribution.

\quad


\section{Introduction}\label{section one}
In this paper, we will give moderate deviation principles from hydrodynamic limits of simple exclusion processes on one dimensional tori starting from  nonequilibrium states, which extend the result given in \cite{Gao2003} about the equilibrium case. We first recall the definition of the simple exclusion process. For each integer $N\geq 1$, let $\mathbb{T}^N=\{0,1,2,\ldots, N-1\}$. The simple exclusion process $\{\hat{\eta}_t\}_{t\geq 0}$ on $\mathbb{T}^N$ is a continuous-time Markov process with state space $\mathbb{X}^N=\{0, 1\}^{\mathbb{T}^N}$ and generator $\hat{\mathcal{L}}_N$ given by
\begin{equation}\label{equ 1.1 generator}
\hat{\mathcal{L}}_{N}f(\eta)=\sum_{x\in \mathbb{T}^N}\left[f(\eta^{x, x+1})-f(\eta)\right]
\end{equation}
for any $\eta\in \mathbb{X}^N$ and $f$ from $\mathbb{X}^N$ to $\mathbb{R}$, where $\eta^{x, x+1}\in \mathbb{X}^N$ is defined as
\[
\eta^{x,x+1}(y)=
\begin{cases}
\eta(y)&\text{~if~}y\neq x, x+1,\\
\eta(x+1)&\text{~if~}y=x,\\
\eta(x) &\text{~if~}y=x+1.
\end{cases}
\]
Note that, throughout this paper operations on $\mathbb{T}^N$ are under the $({\rm mod~}N)$-meaning. For example, $(N-1)+1=0$.

According to Equation \eqref{equ 1.1 generator}, the simple exclusion process describes a traffic flow on $\mathbb{T}^N$. At each $x\in \mathbb{T}^N$, there is at most one particle. All particles perform random walks on $\mathbb{T}^N$. In detail, a particle at $x$ jumps to the neighbor $y=x\pm1$ at rate $1$ when $y$ is not occupied by other particles.

Here we recall invariant measures of simple exclusion processes. For $0<\alpha<1$, let $\nu_\alpha$ be the product measure on $\mathbb{T}^N$ under which $\{\eta(x)\}_{0\leq x\leq N-1}$ are independent and $\nu_\alpha(\eta(x)=1)=\alpha$ for all $0\leq x\leq N-1$, then by Equation \eqref{equ 1.1 generator}, it is easy to check that $\nu_\alpha$ is a reversible distribution and hence an invariant distribution of $\{\hat{\eta}_t\}_{t\geq 1}$. The measure $\nu_\alpha$ is called the global invariant measure of the process. Since the total number of particles is conserved for the simple exclusion process, let
\[
\nu_K(\cdot)=\nu_\alpha\left(\cdot\Big|\sum_{x\in \mathbb{T}^N}\eta(x)=K\right)
\]
for given integer $1\leq K\leq N-1$, then $\nu_K$ is an invariant measure of the process which is called the local invariant measure. Note that $\nu_K$ is independent of the choice of $\alpha$ since the process is irreducible and has finite states conditioned on the total number of particles being $K$.

For other basic properties of simple exclusion processes, readers could see Chapter 8 of \cite{Lig1985} and Part ${\rm \uppercase\expandafter{\romannumeral3}}$ of \cite{Lig1999} for a detailed survey.

Now we recall hydrodynamic limits of simple exclusion processes. We write $\hat{\eta}_t$ as $\hat{\eta}_t^N$ to distinguish different $N$. We use $\eta_t^N$ to denote $\hat{\eta}_{tN^2}^N$. Then, the generator $\mathcal{L}_N$ of $\{\eta_t^N\}_{t\geq 0}$ is given by
$\mathcal{L}_N=N^2\hat{\mathcal{L}}_N$. Let $\mathbb{T}=[0, 1)$ be the one dimensional torus. We denote by $\mu_t^N$ the empirical density field of $\eta_t^N$, i.e.,
\[
\mu_t^N(du)=\frac{1}{N}\sum_{x\in \mathbb{T}^N}\eta_t^N(x)\delta_{\frac{x}{N}}(du),
\]
where $\delta_{a}$ is the Dirac measure concentrated at $a$.  The following hydrodynamic limit theorem is given in Chapter 4 of \cite{kipnis+landim99}.
\begin{proposition}\label{proposition 1.1 hydrodynamic}
(Kipnis and Landim, \cite{kipnis+landim99})
If $\{\eta_0^N(x)\}_{0\leq x\leq N-1}$ are independent and
\[
P(\eta_0^N(x)=1)=\phi\left(\frac{x}{N}\right)
\]
for some $\phi\in C(\mathbb{T})$ and all $N\geq 1, x\in \mathbb{T}^N$, then $\mu_t^N(f)$ converges in probability to $\int_\mathbb{T}\rho_t(u)f(u)du$ as $N\rightarrow+\infty$ for any $t\geq 0$ and $f\in C(\mathbb{T})$, where $\{\rho_t\}_{t\geq 0}$ is the unique weak solution to the heat equation
\begin{equation}\label{equ 1.2 heat equation}
\begin{cases}
&\partial_t\rho(t,u)=\partial_{uu}^2\rho(t,u),\\
&\rho_0=\phi.
\end{cases}
\end{equation}
\end{proposition}
The proof of Proposition \ref{proposition 1.1 hydrodynamic} utilizes Dynkin's martingale formula to show that any weak limit of a subsequence of $\{\mu^N_t\}_{N\geq 1}$ is absolutely continuous with respect to Lebesgue measure and the corresponding R-N derivative is a weak solution to Equation \eqref{equ 1.2 heat equation}.

It is natural to further investigate central limit theorems, large and moderate deviations from hydrodynamic limits given in Proposition \ref{proposition 1.1 hydrodynamic}. The central limit theorem from the hydrodynamic limit is also called the fluctuation. It is shown in Chapter 11 of \cite{kipnis+landim99} that the fluctuation of the simple exclusion process is driven by a generalized Ornstein-Uhlenbeck process introduced in \cite{Holley1978}. The large deviation principle from the hydrodynamic limit of the simple exclusion process is given in \cite{Kipnis1989}, the proof of which utilizes an exponential martingale strategy. A moderate deviation principle is given in \cite{Gao2003} in the case where the initial distribution of the process is $\nu_\alpha$, the proof of which extends the strategy given in \cite{Kipnis1989} and relies heavily on the fact that $\nu_\alpha$ is an invariant measure of the process.

In this paper, we will extend the result given in \cite{Gao2003} to cases where initial distributions of our processes are not invariant. The proof of our main result still utilizes the exponential martingale strategy as that in \cite{Gao2003} but two technical details are improved. We give new approaches to check a replacement lemma and the exponential tightness of the scaled density field of the process, where the invariance of the initial distribution is not need. For mathematical details, see Section \ref{section three}.

\section{Main result}\label{section two}
In this section, we give our main result. For later use, we first introduce some notations and definitions. Let $\mathbb{T}=[0, 1)$ be the one dimensional torus defined as in Section \ref{section one}. Throughout this paper, operators on $\mathbb{T}$ are under $(\text{mod~}1)$-meaning. For example, $0.2-0.3=0.9$. We use $\mathcal{S}$ to denote the dual of $C^\infty(\mathbb{T})$ endowed with the weak topology, i.e., $\nu_n\rightarrow \nu$ in $\mathcal{S}$ when and only when
\[
\lim_{n\rightarrow+\infty}\nu_n(f)=\nu(f)
\]
for any $f\in C^\infty(\mathbb{T})$. Let $\{a_N\}_{N\geq 1}$ be a positive sequence such that
\[
\lim_{N\rightarrow+\infty}\frac{a_N}{N}=\lim_{N\rightarrow+\infty}\frac{\sqrt{N}}{a_N}=0.
\]
Throughout this paper, we adopt the following assumption.

\textbf{Assumption} (A): $\{\eta_0^N(x)\}_{x\in \mathbb{T}^N}$ are independent and
\[
P\left(\eta_0^N(x)=1\right)=\phi\left(\frac{x}{N}\right)
\]
for all $x\in \mathbb{T}^N$, where $\phi\in C(\mathbb{T})$ such that $0<\phi(u)<1$ for all $u\in \mathbb{T}$.

Let $\mathbb{E}$ be the expectation operator. For any $t\geq 0$ and $N\geq 1$, we define the scaled density field $\theta_t^N$ as
\[
\theta_t^N(du)=\frac{1}{a_N}\sum_{x\in \mathbb{T}^N}\left(\eta_t^N(x)-\mathbb{E}\eta_t^N(x)\right)\delta_{\frac{x}{N}}(du).
\]
We can consider $\theta_t^N$ as a random element in $\mathcal{S}$ such that
\[
\theta_t^N(f)=\frac{1}{a_N}\sum_{x\in \mathbb{T}^N}\left(\eta_t^N(x)-\mathbb{E}\eta_t^N(x)\right)f\left(\frac{x}{N}\right).
\]
For given $T>0$, we use $\theta^N$ to denote $\{\theta_t^N\}_{0\leq t\leq T}$. Then $\theta^N$ is a random element in $\mathcal{D}([0, T], \mathcal{S})$, which is the set of c\`{a}dl\`{a}g functions from $[0, T]$ to $\mathcal{S}$ endowed with the Skorokhod topology.

We first give the moderate deviation rate function of the dynamic of the process. For any $W\in \mathcal{D}([0, T], \mathcal{S})$, we define
\begin{align}\label{equ MDP dynamic rate function}
I_{dyn}(W)=&\sup_{F\in C^{1, \infty}([0, T]\times\mathbb{T})}\Bigg\{
W_T(F_T)-W_0(F_0)-\int_0^TW_s((\partial_s+\partial^2_{uu})F_s)ds\notag\\
&-\int_0^T\int_\mathbb{T}\rho_s(u)(1-\rho_s(u))(\partial_uF_s(u))^2dsdu\Bigg\},
\end{align}
where $\{\rho_t\}_{t\geq 0}$ is the unique weak solution to Equation \eqref{equ 1.2 heat equation}. Then, we give the moderate deviation rate function of the intial state of the process. For any $\nu \in \mathcal{S}$, we define
\begin{align}\label{equ MDP intial rate function}
I_{ini}(\nu)=\sup_{f\in C^\infty(\mathbb{T})}\left\{\nu(f)-\frac{1}{2}\int_\mathbb{T}\phi(u)\left(1-\phi(u)\right)f^2(u)du\right\}.
\end{align}
Now we can give our main result.

\begin{theorem}\label{theorem 2.1 main result}
Under Assumption (A),
\begin{equation}\label{equ MDP upper bound}
\limsup_{N\rightarrow+\infty}\frac{N}{a_N^2}\log P\left(\theta^N\in C\right)\leq -\inf_{W\in C}\left(I_{ini}(W_0)+I_{dyn}(W)\right)
\end{equation}
for any closed set $C\subseteq \mathcal{D}([0, T], \mathcal{S})$ and
\begin{equation}\label{equ MDP lower bound}
\liminf_{N\rightarrow+\infty}\frac{N}{a_N^2}\log P\left(\theta^N\in O\right)\geq -\inf_{W\in O}\left(I_{ini}(W_0)+I_{dyn}(W)\right)
\end{equation}
for any open set $O\subseteq \mathcal{D}([0, T], \mathcal{S})$.
\end{theorem}

If $\phi\equiv\alpha$ for some $\alpha\in (0, 1)$, then $\rho_s(u)(1-\rho_s(u))\equiv \alpha(1-\alpha)$ and hence Theorem \ref{theorem 2.1 main result} reduces to Theorem 1.1 of \cite{Gao2003}.

The rest of this paper is devoted to the proof of Theorem \ref{theorem 2.1 main result}. The outline of the proof is as follows. As preliminaries, we first prove a replacement lemma and check the exponential tightness of $\{\theta^N\}_{N\geq 1}$ in the absence of the invariance of the initial distribution. The proof of above replacement lemma utilizes the upper bound of the large deviation principle given in \cite{Kipnis1989}. According to Grownwall's inequality, we reduce the check of the exponential tightness of $\{\theta^N\}_{N\geq 1}$ to that of the logarithm of an exponential martingale and then the required exponential tightness follows from Doob's inequality. For mathematical details, see Section \ref{section three}. With above replacement lemma and exponential tightness, the strategy introduced in \cite{Gao2003} applies to cases discussed in this paper and then Theorem \ref{theorem 2.1 main result} holds. For mathematical details, see Section \ref{section four}.

\section{Exponential tightness and replace lemma}\label{section three}
In this section, we prove following two lemmas.

\begin{lemma}\label{lemma 3.1 exponential tightness}
Under Assumption (A), $\{\theta^N\}_{N\geq 1}$ are exponentially tight.
\end{lemma}

\begin{lemma}\label{lemma 3.2 replacement lemma}
Under Assumption (A), for any $\epsilon>0$ and $G\in C^{1,0}([0, T]\times \mathbb{T})$,
\begin{align}\label{equ 3.1 replacement lemma}
\limsup_{N\rightarrow+\infty}\frac{N}{a_N^2}\log P\Bigg(\Bigg|&\int_0^T\frac{1}{N}\sum_{x\in \mathbb{T}^N}\eta_s^N(x)\eta_s^N(x+1)G_s\left(\frac{x}{N}\right)ds\notag\\
&\text{\quad\quad}-\int_0^T\int_\mathbb{T}\rho_s^2(u)G_s(u)duds\Bigg|>\epsilon\Bigg)=-\infty.
\end{align}

\end{lemma}

Lemmas \ref{lemma 3.1 exponential tightness} and \ref{lemma 3.2 replacement lemma} are analogues of Lemmas 2.1 and 3.2 of \cite{Gao2003} respectively. In \cite{Gao2003}, the simple exclusion process $\eta_t^N$ follows invariant distribution $\nu_\alpha$ at any moment $t$. With this property, Lemma 2.1 of \cite{Gao2003} is proved by utilizing Jensen's inequality and Cauchy-Schwarz inequality. Lemma 3.2 of \cite{Gao2003} is proved by utilizing Garsia-Rademich-Rumsey inequality. For mathematical details, see \cite{Gao2003}. Above approaches do not apply to cases discussed in this paper since initial distributions of simple exclusion processes in this paper are not invariant. Hence we prove above two lemmas in different ways than those in \cite{Gao2003}. By utilizing Grownwall's inequality, we reduce the check of the exponential tightness of $\theta^N$ to that of the logarithm of an exponential martingale and then Lemma \ref{lemma 3.1 exponential tightness} holds according to Doob's inequality. For mathematical details, see Subsection \ref{subsection 3.1}. By utilizing the upper bound of the large deviation principle given in \cite{Kipnis1989}, Lemma \ref{lemma 3.2 replacement lemma} follows from the fact that the minimum of the large deviation rate function given in \cite{Kipnis1989} on a closed set without $\{\rho_t\}_{0\leq t\leq T}$ is strictly positive. For mathematical details, see Subsection \ref{subsection 3.2}.

\subsection{Proof of Lemma \ref{lemma 3.1 exponential tightness}}\label{subsection 3.1}
In this subsection, we prove Lemma \ref{lemma 3.1 exponential tightness}. We first introduce some notations and definitions for later use. For integer $n\geq 0$, we use $e_n(u)$ to denote $\cos(2n\pi u)$. For integer $n\geq 1$, we use $e_{-n}(u)$ to denote $\sin(2n\pi u)$. For $n\in \mathbb{Z}$, $N\geq 1$, $c>0$ and $t\geq 0$, we define
\[
Y_t^{N,n,c}=\exp\left(\frac{a_N^2}{N}\theta_t^N(ce_n)\right).
\]
Furthermore, we define
\[
\mathcal{M}_t^{N,n,c}=\frac{Y_t^{N,n,c}}{Y_0^{N,n,c}}\exp\left(-\int_0^t\frac{(\partial_s+\mathcal{L}_N)Y_s^{N,n,c}}{Y_s^{N,n,c}}ds\right).
\]
By Feynman-Kac formula, $\{\mathcal{M}_t^{N,n,c}\}_{t\geq 0}$ is a martingale with mean $1$.

Now we give the proof of Lemma \ref{lemma 3.1 exponential tightness}

\proof[Proof of Lemma \ref{lemma 3.1 exponential tightness}]

Since ${\rm span}\{e_n:~-\infty<n<+\infty\}$ are dense in $C^\infty(\mathbb{T})$, according to the criterion given in \cite{Puhalskii1994}, to complete this proof we only need to show that
\begin{equation}\label{equ 3.1.1}
\limsup_{M\rightarrow+\infty}\limsup_{N\rightarrow+\infty}\frac{N}{a_N^2}\log P\left(\sup_{0\leq t\leq T}|\theta_t^N(e_n)|>M\right)=-\infty
\end{equation}
and
\begin{equation}\label{equ 3.1.2}
\limsup_{\delta\rightarrow 0}\limsup_{N\rightarrow+\infty}\frac{N}{a_N^2}
\log\sup_{\sigma\in \mathcal{T}}P\left(\sup_{0\leq t\leq \delta}\left|\theta^N_{t+\sigma}(e_n)-\theta^N_\sigma(e_n)\right|>\epsilon\right)=-\infty
\end{equation}
for any $n\in \mathbb{Z}$ and $\epsilon>0$, where $\mathcal{T}$ is the set of stopping times of $\{\eta_t^N\}_{t\geq 0}$ bounded by $T$.

We first check Equation \eqref{equ 3.1.1}. According to Chapman-Kolmogorov equation, we have
\[
\frac{\partial_sY_s^{N,n,1}}{Y_s^{N,n,1}}=-\frac{a_N}{N}\sum_{x\in \mathbb{T}^N}\mathbb{E}\eta_t^N(x)(\Delta^Ne_n(x/N)),
\]
where $\Delta^Ne_n(u)=N^2\left(e_n(u+\frac{1}{N})+e_n(u-\frac{1}{N})-2e_n(u)\right)$. According to the definition of $\mathcal{L}_N$, we have
\[
\frac{\mathcal{L}_NY_s^{N,n,1}}{Y_s^{N,n,1}}=\sum_{x\in \mathbb{T}^N}\left(\exp\left(\frac{a_N}{N}(\eta^N_s(x+1)-\eta_s^N(x))(e_n(x/N)-e_n((x+1)/N))\right)-1\right).
\]
Then, according to the fact that $e^x-1=x+\frac{x^2}{2}+O(x^3)$, we have
\begin{equation*}
\mathcal{M}_t^{N,n,1}=\exp\left(\frac{a_N^2}{N}\left(\theta_t^N(e_n)-\theta_0^N(e_n)-\int_0^t\theta_s^N(\Delta^N e_n)ds-\varepsilon_{1,t}^{N,n}+o(1)\right)\right),
\end{equation*}
where
\[
\varepsilon_{1,t}^{N,n}=\frac{1}{2}\int_0^t\frac{1}{N}\sum_{x\in \mathbb{T}^N}(\eta_s^N(x+1)-\eta_s^N(x))^2(e_n(\frac{x}{n})-e_n(\frac{x+1}{n}))^2ds
\]
and hence $\sup_{t\leq T}|\varepsilon_{1,t}^{N,n}|\leq 2T$. Let $Z_t^{N,n}=\theta_t^N(e_n)-\theta_0^N(e_n)-\int_0^t\theta_s^N(\partial^2_{uu}e_n)ds$, then
\begin{equation}\label{equ 3.1.3}
\mathcal{M}_t^{N,n,1}=\exp\left(\frac{a_N^2}{N}\left(Z_t^{N,n}-\varepsilon_{1,t}^{N,n}+o(1)\right)\right),
\end{equation}
since $|\Delta^Ne_n(u)-\partial^2_{uu}e_n(u)|=O(N^{-1})$. According to the fact that $\partial^2_{uu}e_n=-(2n\pi)^2e_n$ and Grownwall's inequality,
\[
|\theta_t^N(e_n)|\leq \left(|\theta_0^N(e_n)|+\sup_{0\leq t\leq T}Z_t^{N,n}\right)e^{(2n\pi)^2T}
\]
for any $0\leq t\leq T$. According to Assumption (A), it is easy to check that
\[
\limsup_{M\rightarrow+\infty}\limsup_{N\rightarrow+\infty}\frac{N}{a_N^2}\log P\left(|\theta_0^N(e_n)|>M\right)=-\infty.
\]
Hence, to prove Equation \eqref{equ 3.1.1} we only need to show that
\begin{equation}\label{equ 3.1.4}
\limsup_{M\rightarrow+\infty}\limsup_{N\rightarrow+\infty}\frac{N}{a_N^2}\log P\left(\sup_{0\leq t\leq T}Z_t^{N,n}>M\right)=-\infty.
\end{equation}
By Equation \eqref{equ 3.1.3} and Doob's inequality, for sufficiently large $N$,
\begin{align*}
P\left(\sup_{0\leq t\leq T}Z_t^{N,n}>M\right)&\leq P\left(\sup_{0\leq t\leq T}\mathcal{M}_t^{N,n,1}\geq \exp\left(\frac{a_N^2}{N}(M-2T-1)\right)\right)\\
&\leq \exp\left(-\frac{a_N^2}{N}(M-2T-1)\right)
\end{align*}
and hence Equation \eqref{equ 3.1.4} holds. Consequently, Equation \eqref{equ 3.1.1} holds.

Now we check Equation \eqref{equ 3.1.2}. According to the definition of $Z^{N,n}_t$,
\[
\left|\theta_{t+\sigma}^N(e_n)-\theta_\sigma^N(e_n)\right|
\leq \sup_{0\leq s\leq \delta}\left(Z^{N,n}_{\sigma+s}(e_n)-Z^{N,n}_\sigma(e_n)\right)+(2n\pi)^2\delta\sup_{0\leq s\leq T+\delta}|\theta_s^N(e_n)|
\]
for $0\leq t\leq \delta$. Hence, by Equation \eqref{equ 3.1.1},  to prove Equation \eqref{equ 3.1.2} we only need to show that
\begin{equation}\label{equ 3.1.5}
\limsup_{\delta\rightarrow 0}\limsup_{N\rightarrow+\infty}\frac{N}{a_N^2}
\log\sup_{\sigma\in \mathcal{T}}P\left(\sup_{0\leq s\leq \delta}\left(Z^{N,n}_{\sigma+s}(e_n)-Z^{N,n}_\sigma(e_n)\right)>\epsilon\right)=-\infty.
\end{equation}
According to an analysis similar with that leading to Equation \eqref{equ 3.1.3}, we have
\begin{equation}\label{equ 3.1.6}
\frac{\mathcal{M}_{t+\sigma}^{N,n,c}}{\mathcal{M}_\sigma^{N,n,c}}=\exp\left(\frac{a_N^2}{N}\left(c(Z^{N,n}_{\sigma+t}(e_n)-Z^{N,n}_\sigma(e_n))
-c^2\varepsilon_{2,t}^{N,n}+o(1)\right)\right)
\end{equation}
for any $c>0$, where $\varepsilon_{2,t}^{N,n}=\frac{1}{2}\int_\sigma^{\sigma+t}\frac{1}{N}\sum_{x\in \mathbb{T}^N}(\eta_s^N(x+1)-\eta_s^N(x))^2(e_n(\frac{x}{n})-e_n(\frac{x+1}{n}))^2ds$ and hence $|\varepsilon_{2,t}^{N,n}|\leq 2\delta$ for $t\leq \delta$. Therefore, by Doob's inequality,
\begin{align*}
P\left(\sup_{0\leq s\leq \delta}\left(Z^{N,n}_{\sigma+s}(e_n)-Z^{N,n}_\sigma(e_n)\right)>\epsilon\right)&
\leq P\left(\frac{\mathcal{M}_{t+\sigma}^{N,n,c}}{\mathcal{M}_\sigma^{N,n,c}}\geq\exp\left(\frac{a_N^2}{N}\left(c\epsilon-2c^2\delta-\frac{1}{\sqrt{a_N}}\right)\right)
\right)\\
&\leq \exp\left(-\frac{a_N^2}{N}\left(\frac{1}{2}c\epsilon-2c^2\delta\right)\right)
\end{align*}
for sufficiently large $N$ and hence
\[
\limsup_{\delta\rightarrow 0}\limsup_{N\rightarrow+\infty}\frac{N}{a_N^2}
\log\sup_{\sigma\in \mathcal{T}}P\left(\sup_{0\leq s\leq \delta}\left(Z^{N,n}_{\sigma+s}(e_n)-Z^{N,n}_\sigma(e_n)\right)>\epsilon\right)\leq -\frac{1}{2}c\epsilon.
\]
Since $c$ is arbitrary, let $c\rightarrow+\infty$ and then Equation \eqref{equ 3.1.5} holds. Consequently, Equation \eqref{equ 3.1.2} holds and the proof is complete.

\qed

\subsection{Proof of Lemma \ref{lemma 3.2 replacement lemma}}\label{subsection 3.2}
In this subsection, we prove Lemma \ref{lemma 3.2 replacement lemma}. We first recall the upper bound of the large deviation principle given in \cite{Kipnis1989}. We denote by $\mathbb{M}$ the set of measures $\nu$ on $\mathbb{T}$ such that $\nu(\mathbb{T})\leq 1$. Let $\mathbb{M}$ be endowed with the weak topology. For any $\nu\in \mathbb{M}$, let $J_{ini}(\nu)$ be defined as
\[
J_{ini}(\nu)=\sup_{f_1, f_2\in C(\mathbb{T})}\left\{\nu(f_1-f_2)+\int_\mathbb{T}f_2(u)du-\int_\mathbb{T}\log\left(\phi(u)e^{f_1(u)}+(1-\phi(u))e^{f_2(u)}\right)du\right\}.
\]
Let $\mathcal{D}([0, T], \mathbb{M})$ be the set of c\`{a}dl\`{a}g functions from $[0, T]$ to $\mathbb{M}$ endowed with the Skorokhod topology. For any $W\in \mathcal{D}([0, T], \mathbb{M})$. let $J_{dyn}(W)$ be defined as
\begin{align*}
J_{dyn}(W)=&\sup_{F\in C^{1, 2}([0, T]\times\mathbb{T})}\Bigg\{
W_T(F_T)-W_0(F_0)-\int_0^TW_s((\partial_s+\partial^2_{uu})F_s)ds\notag\\
&-\int_0^T\int_\mathbb{T}\frac{dW_s}{du}(u)(1-\frac{dW_s}{du}(u))(\partial_uF_s(u))^2dsdu\Bigg\}
\end{align*}
if $W_s$ is absolutely continuous with respect to Lebesgue measure for all $0\leq s\leq T$ and $J_{dyn}(W)=+\infty$ otherwise. The following proposition is given in \cite{Kipnis1989}.

\begin{proposition}\label{proposition 3.2.1 LDP} (Kipnis, Olla and Varadhan, \cite{Kipnis1989}) Let $\mu_t^N$ be defined as in Section \ref{section one} and $\mu^N=\{\mu_t^N\}_{0\leq t\leq T}$, then
\[
\limsup_{N\rightarrow+\infty}\frac{1}{N}\log P\left(\mu^N\in C\right)\leq -\inf_{W\in C}(J_{dyn}(W)+J_{ini}(W_0))
\]
for any closed set $C\subseteq \mathcal{D}([0, T], \mathbb{M})$.
\end{proposition}

Now we introduce some notations and definitions for later use. For any $t\geq 0$, let $\mu_t(du)=\rho_t(u)du$, where $\{\rho_t\}_{t\geq 0}$ is the unique weak solution to Equation \eqref{equ 1.2 heat equation} defined as in Section \ref{section one}. We further define
\[
\mu=\{\mu_t\}_{0\leq t\leq T}.
\]
For any $f\in C(\mathbb{T})$ and $u\in \mathbb{T}$, let $\tau_uf$ be the element in $C(\mathbb{T})$ such that $\tau_uf(v)=f(v-u)$ for any $v\in \mathbb{T}$. For any $\epsilon>0$ and $f\in C(\mathbb{T})$, we define
\[
\mathcal{C}_{\epsilon, f}=\left\{W\in \mathcal{D}([0, T], \mathbb{M}):~\sup_{u\in \mathbb{T}, \atop 0\leq s\leq T}\left|W_s(\tau_uf)-\mu_s(\tau_uf)\right|\geq \epsilon\right\}.
\]

The following two lemmas are crucial for the proof of Lemma \ref{lemma 3.2 replacement lemma}.

\begin{lemma}\label{lemma 3.2.1}
For any closed set $C\subseteq \mathcal{D}([0, T], \mathbb{M})$, if $C\not\ni \mu$, then
\[
\limsup_{N\rightarrow+\infty}\frac{N}{a_N^2}\log P(\mu^N\in C)=-\infty.
\]
\end{lemma}

\begin{lemma}\label{lemma 3.2.2}
For any $\epsilon>0$ and $f\in C(\mathbb{T})$, $\mathcal{C}_{\epsilon, f}$ is a closed subset of $\mathcal{D}([0, T], \mathbb{M})$.
\end{lemma}

We prove Lemmas \ref{lemma 3.2.1} and \ref{lemma 3.2.2} at the end of this subsection. Now we utilize Lemmas \ref{lemma 3.2.1} and \ref{lemma 3.2.2} to prove Lemma \ref{lemma 3.2 replacement lemma}.

\proof[Proof of Lemma \ref{lemma 3.2 replacement lemma}] The weak solution $\{\rho_t\}_{t\geq 0}$ to Equation \eqref{equ 1.2 heat equation} has a $C([0,+\infty)\times \mathbb{T})$-valued version
\[
\rho_t(u)=\mathbb{E}\phi(u+\sqrt{2}B_t),
\]
where $\{B_t\}_{t\geq 0}$ is the standard Brownian motion starting at $0$.
According to uniform continuities of $G$ and $\{\rho_t\}_{0\leq t\leq T}$, for any $\epsilon>0$, there exists $\delta_1=\delta_1(\epsilon)$ such that
\[
\left|\frac{1}{2\delta N}\sum_{-\delta N\leq j\leq \delta N}G_s\left(\frac{x+j}{N}\right)-G_s\left(\frac{x}{N}\right)\right|\leq \epsilon
\]
and
\[
\left|\frac{1}{2\delta}\int_{-\delta}^\delta \rho_s(u+v)dv-\rho_s(u)\right|\leq \epsilon
\]
for any $\delta\leq \delta_1, x\in \mathbb{T}^N, 0\leq s\leq T, u\in \mathbb{T}$ and sufficiently large $N$. Then, according to the fact that
\begin{align*}
&\frac{1}{N}\sum_{x\in \mathbb{T}^N}\eta_s^N(x)\eta_s^N(x+1)\left(\frac{1}{2\delta N}\sum_{-\delta N\leq j\leq \delta N}G_s\left(\frac{x+j}{N}\right)\right)\\
&=\frac{1}{N}\sum_{x\in \mathbb{T}^N}\left(\frac{1}{2\delta N}\sum_{-\delta N\leq j\leq \delta N}\eta_s^N(x+j)\eta_s^N(x+j+1)\right)G_s\left(\frac{x}{N}\right),
\end{align*}
to prove Lemma \ref{lemma 3.2 replacement lemma} we only need to show that
\begin{align}\label{equ 3.2.5}
\limsup_{\delta\rightarrow0}\limsup_{N\rightarrow+\infty}&\frac{N}{a_N^2}\log P\Bigg(\Bigg|\int_0^T\frac{1}{N}\sum_{x\in \mathbb{T}^N}V_s^N(\eta^N, x, \delta)G_s\left(\frac{x}{N}\right)ds\notag\\
&-\int_0^T\int_\mathbb{T}\left(\frac{1}{2\delta}\int_{-\delta}^{\delta}\rho_s(u+v)dv\right)^2G_s(u)duds\Bigg|>\epsilon\Bigg)=-\infty
\end{align}
for any $\epsilon>0$, where
\[
V_s^N(\eta^N, x, \delta)=\frac{1}{2\delta N}\sum_{-\delta N\leq j\leq \delta N}\eta_s^N(x+j)\eta_s^N(x+j+1).
\]
According to Theorem 2.1 of \cite{Kipnis1989},
\begin{align*}
\limsup_{\delta\rightarrow0}\limsup_{N\rightarrow+\infty}&\frac{1}{N}\log P\Bigg(\Bigg|\int_0^T\frac{1}{N}\sum_{x\in \mathbb{T}^N}V_s^N(\eta^N, x, \delta)G_s\left(\frac{x}{N}\right)ds\notag\\
&-\int_0^T\frac{1}{N}\sum_{x\in \mathbb{T}^N}\left(U_s^N(\eta^N, x, \delta)\right)^2G_s\left(\frac{x}{N}\right)ds\Bigg|>\epsilon\Bigg)=-\infty,
\end{align*}
where
\[
U_s^N(\eta^N, x, \delta)=\frac{1}{2\delta N}\sum_{-\delta N\leq j\leq \delta N}\eta_s^N(x+j)=\mu_s^N(\frac{1}{2\delta}\tau_{\frac{x}{N}}1_{[-\delta, \delta]}).
\]
Hence, to prove Equation \eqref{equ 3.2.5} we only need to show that
\begin{align}\label{equ 3.2.6}
\limsup_{N\rightarrow+\infty}&\frac{N}{a_N^2}\log P\Bigg(\Bigg|\int_0^T\frac{1}{N}\sum_{x\in \mathbb{T}^N}\left(U_s^N(\eta^N, x, \delta)\right)^2G_s\left(\frac{x}{N}\right)ds\notag\\
&-\int_0^T\int_\mathbb{T}\left(\frac{1}{2\delta}\int_{-\delta}^{\delta}\rho_s(u+v)dv\right)^2G_s(u)duds\Bigg|>\epsilon\Bigg)=-\infty
\end{align}
for any $\epsilon, \delta>0$. Note that
\[
\frac{1}{2\delta}\int_{-\delta}^{\delta}\rho_s(u+v)dv=\mu_s(\frac{1}{2\delta}\tau_u1_{[-\delta, \delta]}).
\]
Since
\begin{align*}
&\lim_{N\rightarrow+\infty}\frac{1}{N}\sum_{x\in \mathbb{T}^N}\int_0^T\left(\mu_s\left(\frac{1}{2\delta}\tau_{x/N}1_{[-\delta, \delta]}\right)\right)^2G_s\left(\frac{x}{N}\right)ds\\
&=\int_0^T\int_{\mathbb{T}}\left(\mu_s(\frac{1}{2\delta}\tau_u1_{[-\delta, \delta]})\right)^2G_s(u)dsdu,
\end{align*}
to prove Equation \eqref{equ 3.2.6} we only need to show that
\begin{equation}\label{equ 3.2.7}
\limsup_{N\rightarrow+\infty}\frac{N}{a_N^2}\log P\Bigg(\Bigg|\int_0^T\frac{1}{N}\sum_{x\in \mathbb{T}^N}\mathcal{R}_s^N(\eta^N, x, \delta)G_s\left(\frac{x}{N}\right)ds\Bigg|>\epsilon\Bigg)=-\infty
\end{equation}
for any $\epsilon, \delta>0$, where
\[
\mathcal{R}_s^N(\eta^N, x, \delta)=\left|\left(\mu_s^N(\frac{1}{2\delta}\tau_{x/N}1_{[-\delta,\delta]})\right)^2-\left(\mu_s(\frac{1}{2\delta}\tau_{x/N}1_{[-\delta, \delta]})\right)^2\right|.
\]
For sufficiently large integer $m\geq 1$, let $g_m\in C(\mathbb{T})$ be defined as
\[
g_m(u)=
\begin{cases}
\frac{1}{2\delta} &\text{~if~}|u|\leq \delta,\\
\frac{m}{2\delta}(u+\delta+\frac{1}{m}) & \text{~if~}-\delta-\frac{1}{m}\leq u\leq -\delta,\\
-\frac{m}{2\delta}(x-\delta-\frac{1}{m}) & \text{~if~}\delta\leq u\leq \delta+\frac{1}{m},\\
0 & \text{~else},
\end{cases}
\]
then
\[
\left|\mu_s^N(\frac{1}{2\delta}\tau_{x/N}1_{[-\delta,\delta]})-\mu_s^N(\tau_{x/N}g_m)\right|\leq \frac{1}{\delta m}
\text{~and~}
\left|\mu_s(\frac{1}{2\delta}\tau_{x/N}1_{[-\delta, \delta]})-\mu_s(\tau_{x/N}g_m)\right|\leq \frac{1}{\delta m}.
\]
Hence, to prove Equation \eqref{equ 3.2.7} we only need to show that
\begin{equation}\label{equ 3.2.8}
\limsup_{N\rightarrow+\infty}\frac{N}{a_N^2}\log P\Bigg(\Bigg|\int_0^T\frac{1}{N}\sum_{x\in \mathbb{T}^N}\hat{\mathcal{R}}_s^{N, m}(x)G_s\left(\frac{x}{N}\right)ds\Bigg|>\epsilon\Bigg)=-\infty
\end{equation}
for any $\epsilon>0$ and $m\geq 1$, where
\[
\hat{\mathcal{R}}_s^{N, m}(x)=\left|\left(\mu_s^N(\tau_{x/N}g_m)\right)^2-\left(\mu_s(\tau_{x/N}g_m)\right)^2\right|.
\]
Since $\mu^N_s(\mathbb{T}), \mu_s(\mathbb{T})\leq 1$, to prove Equation \eqref{equ 3.2.8} we only need to show that
\begin{equation}\label{equ 3.2.9}
\limsup_{N\rightarrow+\infty}\frac{N}{a_N^2}\log P\left(\mu^N\in \mathcal{C}_{\epsilon, g_m}\right)=-\infty
\end{equation}
for any $\epsilon>0$ and $m\geq 1$. Since $\mu\not\in \mathcal{C}_{\epsilon, g_m}$, Equation \eqref{equ 3.2.9} follows from Lemmas \ref{lemma 3.2.1} and \ref{lemma 3.2.2}. Consequently, the proof is complete.

\qed

At last, we prove Lemmas \ref{lemma 3.2.1} and \ref{lemma 3.2.2}.

\proof[Proof of Lemma \ref{lemma 3.2.1}]

Since $\lim_{N\rightarrow+\infty}\frac{N}{a_N}=+\infty$, to prove Lemma \ref{lemma 3.2.1}, we only need to show that
\begin{equation}\label{equ 3.2.1}
\limsup_{N\rightarrow+\infty}\frac{1}{N}\log P(\mu^N\in C)<0.
\end{equation}
According to Proposition \ref{proposition 3.2.1 LDP}, to prove Equation \eqref{equ 3.2.1}, we only need to show that
\begin{equation}\label{equ 3.2.2}
\inf_{W\in C}(J_{ini}(W_0)+J_{dyn}(W))>0.
\end{equation}
It is shown in \cite{Kipnis1989} that $J_{ini}(\cdot_0)+J_{dyn}(\cdot)$ is a good rate function. Hence, to prove Equation \eqref{equ 3.2.2}, we only need to show that $J_{ini}(W_0)=J_{dyn}(W)=0$ implies that $W=\mu$. For $W$ making $J_{ini}(W_0)=J_{dyn}(W)=0$, we define
\begin{align*}
l_1(\epsilon, F)=&W_T(\epsilon F_T)-W_0(\epsilon F_0)-\int_0^TW_s((\partial_s+\partial^2_{uu})(\epsilon F_s))ds\\
&-\int_0^T\int_\mathbb{T}\frac{dW_s}{du}(u)(1-\frac{dW_s}{du}(u))(\partial_u(\epsilon F_s)(u))^2dsdu
\end{align*}
for any $\epsilon\in \mathbb{R}$ and $F\in C^{1, 2}([0, T]\times\mathbb{T})$ and
\begin{align*}
&l_2(\epsilon_1, \epsilon_2, f_1, f_2)\\
&=W_0(\epsilon_1f_1-\epsilon_2f_2)+\int_\mathbb{T}\epsilon_2f_2(u)du-\int_\mathbb{T}\log\left(\phi(u)e^{\epsilon_1f_1(u)}+(1-\phi(u))e^{\epsilon_2f_2(u)}\right)du
\end{align*}
for any $\epsilon_1, \epsilon_2\in \mathbb{R}$ and $f_1, f_2\in C(\mathbb{T})$. Since $J_{ini}(W_0)=0$, we have
\[
\partial_{\epsilon_1}l_2(0,0, f_1, f_2)=\partial_{\epsilon_2}l_2(0,0,f_1, f_2)=0
\]
and hence $W_0(f_1)=\int_\mathbb{T}\phi(u)f_1(u)du$ for any $f_1\in C(\mathbb{T})$. Consequently, $W_0(du)=\phi(u)du=\mu_0(du)$. Similarly, $J_{dyn}(W)=0$ implies that
\[
\partial_\epsilon l_1(0, F)=0.
\]
We choose $F$ with the form $F(s,u)=h(s)f(u)$ for some $h\in C^1([0, T])$ and $f\in C^\infty(\mathbb{T})$, then we have
\[
h(T)W_T(f)-h(0)W_0(f)-\int_0^T\partial_sh(s)W_s(f)ds=\int_0^Th(s)W_s(\partial^2_{uu}f)ds.
\]
Since $h$ is arbitrary, $\{W_t(f)\}_{0\leq t\leq T}$ is absolutely continuous and
\begin{equation}\label{equ 3.2.3}
\partial_tW_t(f)=W_t(\partial^2_{uu}f)
\end{equation}
for any $f\in C(\mathbb{T})$. Let $e_n(u)=\cos(2n\pi u)$ for $n\geq 0$ and $e_{-n}(u)=\sin(2n\pi u)$ for $n\geq 1$ defined as in Subsection \ref{subsection 3.1}, then Equation \eqref{equ 3.2.3} and the fact $W_0=\mu_0$ implies that
\[
W_t(e_m)=\mu_0(e_m)e^{-(2m\pi)^2t}
\]
for any integer $m$. Since the span of $\{e_m\}_{-\infty<m<+\infty}$ is dense in $C(\mathbb{T})$, the solution to Equation \eqref{equ 3.2.3} with initial condition $W_0=\mu_0$ is unique. Since $\mu$ is also a solution to Equation \eqref{equ 3.2.3}, we have $W=\mu$ and the proof is complete.

\qed

\proof[Proof of Lemma \ref{lemma 3.2.2}]

Assuming that $W^n\in \mathcal{C}_{\epsilon, f}$ for $n\geq 1$ and $W^n\rightarrow W$ in $\mathcal{D}([0, T], \mathbb{M})$, then we only need to show that $W\in \mathcal{C}_{\epsilon, f}$. Since $W^n\in \mathcal{C}_{\epsilon, f}$ for $n\geq 1$, there exist a sequence $\{v^n\}_{n\geq 1}$ in $\mathbb{T}$ such that
\[
\liminf_{n\rightarrow+\infty}\sup_{0\leq s\leq T}\left|W^n_s(\tau_{v^n}f)-\mu_s(\tau_{v^n}f)\right|\geq \epsilon.
\]
Since $\mathbb{T}$ is compact, there exists $v\in \mathbb{T}$ such that $v$ is the limit of a subsequence of $\{v^n\}_{n\geq 1}$ in $\mathbb{T}$. For simplicity, we still write this subsequence as $\{v^n\}_{n\geq 1}$. Since $W^n\rightarrow W$ in $\mathcal{D}([0, T], \mathbb{M})$ as $n\rightarrow+\infty$, there exist a sequence of increasing continuous functions $\{\varphi^n\}_{n\geq 1}$ from $[0, T]$ to $[0, T]$ such that $\varphi^n(0)=0, \varphi^n(T)=T$ for all $n$, $\varphi^n(s)\rightarrow s$ uniformly for $s\in [0, T]$ and $W^n_{\varphi^n(s)}\rightarrow W_s$ uniformly for $s\in [0, T]$. According to triangle inequality, for any $n\geq 1$,
\begin{align}\label{equ 3.2.4}
\sup_{0\leq s\leq T}\left|W_s(\tau_vf)-\mu_s(\tau_vf)\right|\leq & \sup_{0\leq s\leq T}\left|W_s(\tau_vf)-W^n_{\varphi^n(s)}(\tau_vf)\right| \notag\\
&+\sup_{0\leq s\leq T}\left|W^n_{\varphi^n(s)}(\tau_vf)-W^n_{\varphi^n(s)}(\tau_{v^n}f)\right|\notag\\
&+\sup_{0\leq s\leq T}\left|W^n_{\varphi^n(s)}(\tau_{v^n}f)-\mu_{\varphi^n(s)}(\tau_{v^n}f)\right|\notag\\
&+\sup_{0\leq s\leq T}\left|\mu_{\varphi^n(s)}(\tau_{v^n}f)-\mu_s(\tau_vf)\right|.
\end{align}
According to the definition of $\mu$, $\mu_t$ is continuous in $t$. According to the uniform continuity of $f$, $\tau_{v^n}f(u)\rightarrow \tau_vf(u)$ uniformly for $u\in \mathbb{T}$. Then, let $n\rightarrow+\infty$ in Equation \eqref{equ 3.2.4}, we have
\[
\sup_{0\leq s\leq T}\left|W_s(\tau_vf)-\mu_s(\tau_vf)\right|\geq \epsilon
\]
and hence $W\in \mathcal{C}_{\epsilon, f}$, which completes the proof.

\qed

\section{Proof of Theorem \ref{theorem 2.1 main result}}\label{section four}
In this section, we prove our main result Theorem \ref{theorem 2.1 main result}. With Lemmas \ref{lemma 3.1 exponential tightness} and \ref{lemma 3.2 replacement lemma}, the strategy introduced in \cite{Gao2003} applies to cases discussed in this paper. So we only give outlines of proofs of Equations \eqref{equ MDP upper bound} and \eqref{equ MDP lower bound} to avoid repeating many similar details with those in \cite{Gao2003}. For later use, we define
\[
Y_t^N(F)=\exp\left(\frac{a_N^2}{N}\theta_t^N(F_t)\right)
\]
for any $0\leq t\leq T$ and $F\in C^{1,+\infty}([0, T]\times\mathbb{T})$. Furthermore, we define
\[
\mathcal{M}_t^N(F)=\frac{Y_t^N(F)}{Y_0^N(F)}\exp\left(-\int_0^t\frac{(\partial_s+\mathcal{L}_N)Y_s^N(F)}{Y_s^N(F)}ds\right)
\]
for $0\leq t\leq T$, then $\{\mathcal{M}_t^N(F)\}_{0\leq t\leq T}$ is a martingale according to Feynman-Kac formula. For any $g\in C(\mathbb{T})$, let $P^N_g$ be the probability measure of our simple exclusion process with initial condition where $\{\eta_0^N(x)\}_{x\in \mathbb{T}^N}$ are independent and
\[
P\left(\eta_0^N(x)=1\right)=\phi(x/N)+\frac{a_N}{N}g(x/N)
\]
for all $x\in \mathbb{T}^N$. Then, for any $F\in C^{1,\infty}([0, T]\times\mathbb{T})$, we define $\hat{P}_g^{N,F}$ as the probability measure such that
\[
\frac{d\hat{P}_g^{N, F}}{dP_g^N}=\mathcal{M}_T^N(F).
\]

Now we prove Equation \eqref{equ MDP upper bound}.
\proof[Proof of Equation \eqref{equ MDP upper bound}]

According to Lagrange's mean value theorem and an analysis similar with that leading to Equation \eqref{equ 3.1.3}, we have
\begin{equation}\label{equ 4.1}
\mathcal{M}_T^N(F)=\exp\left(\frac{a_N^2}{N}\left(Z_T^N(F)-\varepsilon_{3,T}^N(F)+o(1)\right)\right),
\end{equation}
where
\[
Z_T^N(F)=\theta_T^N(F_T)-\theta_0(F_0)-\int_0^T\theta_s^N((\partial_s+\partial^2_{uu})F_s)ds
\]
and
\begin{align*}
\varepsilon_{3,T}^N(F)&=\frac{1}{2}\int_0^T\frac{1}{N}\sum_{x\in \mathbb{T}^N}\left(\partial_uF_s(x/N)\right)^2\left(\eta_s^N(x)+\eta_s^N(x+1)-2\eta_s^N(x)\eta_s^N(x+1)\right)ds\\
&=\int_0^T\mu_s^N((\partial_uF_s)^2)-\frac{1}{N}\sum_{x\in \mathbb{T}^N}(\partial_uF_s(x/N))^2\eta_s^N(x)\eta_s^N(x+1)ds+o(1).
\end{align*}
According to Lemmas \ref{lemma 3.1 exponential tightness} and \ref{lemma 3.2.1}, we have
\begin{align*}
\varepsilon_{3,T}^N(F)&=\int_0^T\mu_s((\partial_uF_s)^2)ds-\int_0^T\int_\mathbb{T}\rho_s^2(u)(\partial_uF_s(u))^2duds+\varepsilon_{4,T}^N(F)\\
&=\int_0^T\rho_s(u)(1-\rho_s(u))(\partial_uF_s(u))^2ds+\varepsilon_{4,T}^N(F),
\end{align*}
where
\begin{equation}\label{equ 4.6}
\limsup_{N\rightarrow+\infty}\frac{N}{a_N^2}\log P\left(|\varepsilon_{4,T}^N(F)|>\epsilon\right)=-\infty
\end{equation}
for any $\epsilon>0$. As a result, by Equation \eqref{equ 4.1}, we have
\begin{equation}\label{equ 4.2}
\mathcal{M}_T^N(F)=\exp\left(\frac{a_N^2}{N}\left(\mathcal{I}(\theta^N, F)+\varepsilon_{4,T}^N(F)+o(1)\right)\right),
\end{equation}
where
\begin{align*}
\mathcal{I}(W,F)=&W_T(F_T)-W_0(F_0)-\int_0^TW_s((\partial_s+\partial^2_{uu})F_s)ds\\
&-\int_0^T\int_\mathbb{T}\rho_s(u)(1-\rho_s(u))(\partial_uF_s(u))^2dsdu
\end{align*}
for any $W\in \mathcal{D}([0, T], \mathcal{S})$. According to Assumption (A), it is easy to check that
\begin{equation}\label{equ 4.3}
\mathbb{E}\exp\left(\frac{a_N^2}{N}\theta_0^N(f)\right)=\exp\left(\frac{a_N^2}{2N}\left(\int_\mathbb{T}f^2(u)\phi(u)(1-\phi(u))du+o(1)\right)\right)
\end{equation}
for any $f\in C^\infty(\mathbb{T})$. By utilizing Markov's inequality and the minimax theorem given in \cite{Sion1958}, Equation \eqref{equ MDP upper bound} holds for all compact $C\subseteq \mathcal{D}([0, T], \mathcal{S})$ according to Equations \eqref{equ 4.2} and \eqref{equ 4.3}. To show that Equation \eqref{equ MDP upper bound} holds for all closed $C$, we only need to show that $\{\theta^N\}_{N\geq 1}$ are exponential tight and hence the proof is complete according to Lemma \ref{lemma 3.1 exponential tightness}.

\qed

To prove Equation \eqref{equ MDP lower bound}, we need following two lemmas.

\begin{lemma}\label{lemma 4.1}
If $W$ makes $I_{ini}(W_0)+I_{dyn}(W)<+\infty$, then there exist $g,F$ such that $W_t(du)=\rho^{F,g}(t,u)du$ for all $0\leq t\leq T$, where
\begin{equation*}
\begin{cases}
&\frac{d}{dt}\rho^{F,g}(t,u)=\partial^2_{uu}\rho^{F,g}(t,u)-2\frac{\partial}{\partial u}\left(\rho_t(u)(1-\rho_t(u))\partial_uF_t(u)\right),\\
&\rho^{F,g}_0=g.
\end{cases}
\end{equation*}
Furthermore, $I_{ini}(W_0)=\frac{1}{2}\int_\mathbb{T}\frac{g^2(u)}{\phi(u)(1-\phi(u))}du$ and
\[
I_{dyn}(W)=\mathcal{I}(W, F)=\int_0^T\int_\mathbb{T}\rho_s(u)(1-\rho_s(u))(\partial_uF_s(u))^2dsdu.
\]
\end{lemma}

\begin{lemma}\label{lemma 4.2}
As $N\rightarrow+\infty$, $\theta^N$ converges in $\hat{P}^{N,F}_g$-probability to $\{\rho^{F,g}(t,u)(du)\}_{0\leq t\leq T}$.
\end{lemma}

Lemmas \ref{lemma 4.1} and \ref{lemma 4.2} are analogues of Lemma 5.1 and Theorem 4.1 of \cite{Gao2003} respectively. With Lemmas \ref{lemma 3.1 exponential tightness} and \ref{lemma 3.2 replacement lemma}, analyses given in proofs of Lemma 5.1 and Theorem 4.1 of \cite{Gao2003} apply to Lemmas \ref{lemma 4.1} and \ref{lemma 4.2} respectively. Hence we omit proofs of Lemmas \ref{lemma 4.1} and \ref{lemma 4.2} here. At last, we prove Equation \eqref{equ MDP lower bound}.

\proof[Proof of Equation \eqref{equ MDP lower bound}]

We only deal with the case where $\inf_{W\in O}\left(I_{ini}(W_0)+I_{dyn}(W)\right)<+\infty$. Otherwise, Equation \eqref{equ MDP lower bound} is trivial. For any $\epsilon>0$, according to Lemma \ref{lemma 4.1}, there exists $W^\epsilon\in O$ such that
\[
I_{ini}(W^\epsilon_0)+I_{dyn}(W^\epsilon)<\inf_{W\in O}\left(I_{ini}(W_0)+I_{dyn}(W)\right)+\epsilon
\]
and $W^\epsilon_t(du)=\rho^{F^\epsilon, g^\epsilon}(t,u)du$ for some $F^\epsilon, g^\epsilon$ and $0\leq t\leq T$. According to Lemma \ref{lemma 4.2}, $\theta^N$ converges in $\hat{P}^{N,F^\epsilon}_{g^\epsilon}$-probability to $\{\rho^{F^\epsilon, g^\epsilon}(t,u)du\}_{0\leq t\leq T}=W^\epsilon$. Then, we have
\begin{equation}\label{equ 4.4}
\lim_{N\rightarrow+\infty}\hat{P}^{N,F^\epsilon}_{g^\epsilon}\left(\theta^N\in O\right)=1
\end{equation}
and
\begin{equation}\label{equ 4.5}
\lim_{N\rightarrow+\infty}\hat{P}^{N, F^\epsilon}_{g^\epsilon}\left(|\mathcal{I}(\theta^N, F^\epsilon)-\mathcal{I}(W^\epsilon, F^\epsilon)|\leq \epsilon\right)=1.
\end{equation}
According to Equation \eqref{equ 4.6}, Assumption (A) and Cauchy-Schwarz inequality, it is easy to check that
\[
\lim_{N\rightarrow+\infty}\hat{P}^{N, F^\epsilon}_{g^\epsilon}\left(|\varepsilon_{4,T}^N(F^\epsilon)|\leq\epsilon\right)=1.
\]
Then, by Equations \eqref{equ 4.2}, \eqref{equ 4.4}, \eqref{equ 4.5} and Lemma \ref{lemma 4.1}, we have
\begin{equation}\label{equ 4.7}
\lim_{N\rightarrow+\infty}\hat{P}^{N,F^\epsilon}_{g^\epsilon}\left(\theta^N\in O, \mathcal{M}_T^N(F^\epsilon)\leq e^{\frac{a_N^2}{N}\left(I_{dyn}(W^\epsilon)+3\epsilon\right)}\right)=1.
\end{equation}
Under Assumption (A), it is easy to check that $\frac{N}{a_N^2}\log\frac{dP}{dP^N_{g^\epsilon}}$ converges in $\hat{P}^{N, F^\epsilon}_{g^\epsilon}$-probability to
\[
-\frac{1}{2}\int_\mathbb{T}\frac{(g^\epsilon(u))^2}{\phi(u)(1-\phi(u))}du=-I_{ini}(W^\epsilon_0)
\]
as $N\rightarrow+\infty$. Hence, by Equation \eqref{equ 4.7},
\begin{align}\label{equ 4.8}
\lim_{N\rightarrow+\infty}\hat{P}^{N,F^\epsilon}_{g^\epsilon}\Big(\theta^N\in O, &\mathcal{M}_T^N(F^\epsilon)\leq e^{\frac{a_N^2}{N}\left(I_{dyn}(W^\epsilon)+3\epsilon\right)},\notag\\
&\frac{dP}{dP^N_{g^\epsilon}}\geq \exp\left(-\frac{a_N^2}{N}(I_{ini}(W^\epsilon_0)+\epsilon)\right)\Big)=1.
\end{align}
According to the definition of $\hat{P}^{N, F^\epsilon}_{g^\epsilon}$, by Equation \eqref{equ 4.8},
\begin{align*}
P\left(\theta^N\in O\right)&=\mathbb{E}_{\hat{P}^{N, F^\epsilon}_{g^\epsilon}}\left(\frac{dP}{dP^N_{g^\epsilon}}\left(\mathcal{M}_T^N(F^\epsilon)\right)^{-1}1_{\{\theta^N\in O\}}\right)\\
&\geq \exp\left(-\frac{a_N^2}{N}\left(I_{ini}(W^\epsilon_0)+I_{dyn}(W^\epsilon)+4\epsilon\right)\right)(1+o(1)).
\end{align*}
Hence,
\begin{align*}
\liminf_{N\rightarrow+\infty}\frac{N}{a_N^2}\log P\left(\theta^N\in O\right)&\geq -(I_{ini}(W^\epsilon_0)+I_{dyn}(W^\epsilon))-4\epsilon\\
&=-\inf_{W\in O}\left(I_{ini}(W_0)+I_{dyn}(W)\right)-5\epsilon.
\end{align*}
Since $\epsilon$ is arbitrary, let $\epsilon\rightarrow 0$ and the proof is complete.

\qed

\quad

\textbf{Acknowledgments.}
The author is grateful to the financial
support from the Fundamental Research Funds for the Central Universities with grant number 2022JBMC039.

{}

\begin{thebibliography}{}
\bibitem{Gao2003}Gao, FQ. and Quastel, J. (2003). Moderate deviations from the hydrodynamic limit of the symmetric exclusion process. \emph{Science in China (Series A)} \textbf{5}, 577-592.
\bibitem{Holley1978}Holley, R. A. and Stroock, D. W. (1978). Generalized Ornstein-Uhlenbeck processes and infinite
particle branching brownian motions. \emph{Publications of the Research Institute for Mathematical
Sciences} \textbf{14}, 741-788.
\bibitem{kipnis+landim99} Kipnis, C. and Landim, C. (1999) {\it Scaling limits of interacting particle systems.} Springer-Verlag, Berlin.
\bibitem{Kipnis1989} Kipnis, C., Olla, S. and Varadhan, S. R. S. (1989). Hydrodynamics and large deviation for simple exclusion processes. \emph{Communications on Pure $\&$ Applied Mathematics} \textbf{42}, 115-137.
\bibitem{Lig1985}Liggett, T. M. (1985). \emph{Interacting Particle Systems.} Springer, New York.
\bibitem{Lig1999}Liggett, T. M. (1999). \emph{Stochastic interacting systems: contact, voter and exclusion processes.}
Springer, New York.
\bibitem{Puhalskii1994}Puhalskii, A. (1994). The method of stochastic exponentials for large deviations. \emph{Stochastic Processes and their Applications} \textbf{54}, 45-70.
\bibitem{Sion1958}Sion, M. (1958). On general minimax theorems. \emph{Pacific Journal of Mathematics} \textbf{8}, 171-176.
\end{thebibliography}
\end{document}